\documentclass[11pt]{article}
\usepackage{amssymb}
\usepackage{graphicx}
\usepackage{amsmath}

\newtheorem{theorem}{Theorem}

\newtheorem{ob}[theorem]{Observation}

\newtheorem{corollary}[theorem]{Corollary}

\newtheorem{lemma}[theorem]{Lemma}

\newenvironment{proof}[1][Proof]
{\textbf{#1.} }{\  \rule{0.5em}{0.5em}}
\input{tcilatex}
\begin{document}

\title{\textbf{Trees with unique minimum glolal offensive alliance sets}}
\author{\textbf{Mohamed Bouzefrane} \\
{\small Faculty of Technology, University of M\'{e}d\'{e}a, Algeria}\\
{\small email: mohamedbouzefrane@gmail.com} \and  \textbf{Isma Bouchemakh} \\
{\small Faculty of Mathematics, Laboratory L'IFORCE,}\\
{\small University of Sciences and Technology Houari Boumediene (USTHB), }\\
{\small B.P. 32 El-Alia, Bab-Ezzouar, 16111 Algiers, Algeria}\\
{\small emails: \ isma\_bouchemakh2001@yahoo.fr}, {\small %
ibouchemakh@usthb.dz} \and \textbf{Mohamed Zamime} \\
{\small Faculty of Technology, University of M\'{e}d\'{e}a, Algeria}\\
{\small email: \ zamimemohamed@yahoo.com} \and \textbf{Noureddine Ikhlef-Eschouf} \\
{\small Faculty of Sciences, Department of Mathematics and Computer Science,}%
\\
{\small University of M\'{e}d\'{e}a, Algeria}\\
{\small email: \ nour\_echouf@yahoo.fr}\ }
\date{}
\maketitle

\begin{abstract}
Let $G=$ $\left( V,E\right) $ be a simple graph.\ A non-empty set $S \subseteq V$ is
called a global offensive alliance if $S$ is a dominating set and for every
vertex $v$ in $V-S$, at least half of the vertices from the closed neighborhood
of $v$ are in $S$. The global offensive alliance number is the minimum
cardinality of a global offensive alliance in $G$. In this paper, we give a
constructive characterization of trees having a unique minimum global
offensive alliance.

\medskip

\end{abstract}


\section{Introduction}

Throughout this paper, $G=(V,E)$ denotes a simple graph with vertex-set $%
V=V(G)$ and edge-set $E=E(G)$. Let $G$ and $H$ be two graphs with two
disjoint vertex sets. Their \emph{disjoint union} is denoted by $G\cup H$,
the disjoint union of $k$ copies of $G$ is denoted by $kG$ and the disjoint
union of a family of graphs $G_{1},G_{2},\ldots,G_{k}$ is denoted by $\cup
_{i=1}^{k}G_{i}.$ For every vertex $v\in V(G),$ the $\emph{open}$ $\emph{%
neighborhood}$ $N_{G}(v)$ is the set $\left \{ u\in V(G)\mid uv\in
E(G)\right \} $ and the \emph{closed neighborhood} of $v$ is the set $N_{G}%
\left[ v\right] =N(v)\cup \left \{ v\right \} .$ The \emph{degree} of a
vertex $v\in V(G),$ denoted $d_{G}\left( v\right) $, is the size of its open
neighborhood. A vertex of degree one is called a \emph{leaf} and its
neighbor is called a \emph{support vertex}. If $v$\ is a support vertex of a
tree $T$, then $L_{T}(v)$\ will denote the set of the leaves attached at $v$%
. Let $L(T)$ and $S(T)$ denote the set of leaves and support vertices,
respectively, in $T,$ and let $\left \vert L(T)\right \vert =l\left(
T\right) $. As usual, the \emph{path} of order $n$ is denoted by $P_{n},$
and the \emph{star} of order $n$ by $K_{1,n-1}.$ A \emph{double star} $%
S_{p,q}$ is obtained by attaching $p$ leaves at an endvertex of a path $%
P_{2} $ and $q$ leaves at the second one. A \emph{subdivision} of an edge $%
uv $ is obtained by introducing a new vertex $w$ and replacing the edge $uv$
with the edges $uw$ and $wv.$ A subdivided star denoted by $SS_{k}$ is a
star $K_{1,k}$ where each edge is subdivided exactly once. A \emph{wounded
spider} is a tree obtained from $K_{1,r},$ where $r\geq 1$, by subdividing at most
$r-1$ of its edges. For a vertex $v,$ let $C(v)$ and $D(v)$ denote the set
of \emph{children} and \emph{descendants}, respectively, of $v$ in a rooted
tree $T$, and let $D[v]=D(v)\cup \{v\}.$ The \emph{maximal subtree} at $v$
is the subtree of $T$ induced by $D[v]$, and is denoted by $T_{v}.$

A \emph{dominating set} of a graph $G$ is a set $D$ of vertices such that
every vertex in $V-D$ is adjacent to some vertex in $D$. \emph{The
domination number} of $G$, denoted by $\gamma \left( G\right)$, is the
minimum cardinality of a dominating set of $G$. The concept of domination in
graphs, with its many variations, is now well studied in graph theory. For
more details, see the books of Haynes, Hedetniemi, and Slater \cite%
{hhs1,hhs2}.

Among the many variations of domination, we mention the concept of alliances
in graphs that has been studied in recent years. Several types of alliances
in graphs are introduced in \cite{hhk1}, including the offensive alliance
that we study here. A dominating set $D$ with the property that for every
vertex $v$ not in $D,$
\begin{equation}
\left \vert N_{G}\left[ v\right] \cap D\right \vert \geq \left \vert N_{G}%
\left[ v\right] -D\right \vert  \label{def}
\end{equation}%
is called \emph{global offensive alliance set of G} and abbreviated \emph{%
GOA-set} of $G$. The \emph{global offensive alliance number} $\gamma
_{o}\left( G\right) $\ is the minimum cardinality among all GOA-sets of $G$.
A GOA-set of $G$ of cardinality $\gamma _{o}\left( G\right) $ is called $%
\gamma _{o}$-set of $G$, or $\gamma _{o}\left( G\right) $-set. Several works
have been carried out on global offensive alliances in graphs (see, for
example, \cite{bc,cv1}, and elsewhere).

Graphs with unique minimum $\mu $-set, where $\mu $ is a some graph
parameter, is another concept to which much attention was given during the
last two decades. For example, graphs with unique minimum $\gamma $-set were
first studied by Gunther et al. in \cite{ghmr1}. Later this problem was
studied for various classes of graphs including block graphs \cite{f1},
cactus graphs \cite{fv1}, some cartesian product graphs \cite{h1} and some
repeated cartesian products \cite{h2}. Several works on uniqueness related
to other graph parameters have been widely studied, such as
locating-domination number \cite{bcl}, paired-domination number \cite{ch1},
double domination number \cite{ch2}, roman domination number \cite{cr1} and
total domination number \cite{hh1}. Further work on this topic can be found
in \cite{frv1,fv2,fvz1,fs1,h3,hs1,stv1,t1}

The aim of this paper is to characterize all trees having unique minimum
global offensive alliance set. We denote such trees as \emph{UGOA-trees}.

\section{Preliminaries results}

We give in this section the following observations. Some results are
straightforward and so their proofs are omitted.

\begin{ob}
\label{observation 1}Let $T$ be a tree of order at least three and $u\in
S(T).$ Then,

\begin{itemize}
\item[$(i)$] there is a $\gamma _{o}\left( T\right) $-set that contains all
support vertices of $T$,

\item[$(ii)$] if $D$ is a unique $\gamma _{o}\left( T\right) $-set, then $D$
contains all support vertices but no leaf,

\item[$(iii)$] if $l_{T}(u)\geq 2,$ then $u$ belongs to any $\gamma _{o}$%
-set($T$).
\end{itemize}
\end{ob}

\noindent\begin{proof}
$(i)$ and $(ii)$ are obvious. If $(iii)$ is not satisfied, then  all leaves attached at $u$ would be
contained in $D$, which is a contradiction with the minimality of $D$.
\end{proof}

\begin{ob}
\label{observation 2}Let $T$ be a tree obtained from a nontrivial tree $%
T^{\prime }$ by joining a new vertex $v$ at a support vertex $u$ of $%
T^{\prime }.$ Let $D$ and $D^{\prime }$ be $\gamma _{o}\left( T\right) $%
-sets of $T$ and $T^{\prime },$ respectively. Then,

\begin{itemize}
\item[$(i)$] $\left \vert D^{\prime }\right \vert =\left \vert D\right \vert,$

\item[$ii)$] $D\cap V(T^{\prime })$ is a $\gamma _{o}\left( T^{\prime
}\right) $-set,

\item[$(iii)$] if $T$ is a UGOA-tree such that $u$ is in any $\gamma
_{o}\left( T^{\prime }\right) $-set$,$ then $T^{\prime }$ is a UGOA-tree.
\end{itemize}
\end{ob}

\noindent\begin{proof}
According to Observation \ref{observation 1} $(iii),$ $u$ must be in $D$
since $l_{T}(u)\geq 2.$\vspace{2mm}\\
$i)$  $D$ is clearly a GOA-set of $T^{\prime },$ and then $%
\left
\vert D^{\prime }\right \vert \leq \left \vert D\right \vert .$ By
Observation \ref{observation 1} $(i)$, we can assume that $u\in D^{\prime }$. Hence, $D^{\prime }$ can be extended to a GOA-set of $T$, which leads to $%
\left \vert D\right \vert \leq \left \vert D^{\prime }\right \vert .$\ Thus
equality holds.\vspace{2mm}\\
$ii)$ Since $D\cap V(T^{\prime })=D$ is a GOA-set of $T^{\prime }$
with cardinality $\left \vert D\right \vert =\left \vert D^{\prime
}\right\vert$,  we deduce that $D\cap V(T^{\prime })$ is a $\gamma _{o}\left( T^{\prime
}\right) $-set.\vspace{2mm}\\
$iii)$ Item $(i)$ together with the fact that $u$ belongs to any $\gamma
_{o}\left( T^{\prime }\right) $ imply that $D^{\prime }$ can be extended to
a $\gamma _{o}\left( T\right) $-set. Therefore, the uniqueness of $D$ as a $%
\gamma _{o}\left( T\right) $-set leads to $D^{\prime }=D,$ which means that $%
D^{\prime }$ is the unique $\gamma _{o}\left( T^{\prime }\right) .$
\end{proof}

\begin{ob}
\label{observation 3} Let $T$ be a tree obtained from a nontrivial tree $%
T^{\prime }$ different from $P_{2}$ by joining the center vertex $y$ of the
path $P_{3}=x$-$y$-$z$ at a support vertex $v$ of $T^{\prime }$. Let $D$ and
$D^{\prime }$ be $\gamma _{o}\left( T\right) $-sets of $T$ and $T^{\prime },$
respectively such that each of them contains all support vertices. Then,

\begin{itemize}
\item[$(i)$] $\left \vert D^{\prime }\right \vert =\left \vert D\right \vert
-1$,

\item[($ii)$] $D\cap V(T^{\prime })$ is a $\gamma _{o}\left( T^{\prime
}\right) $-set,

\item[$(iii)$] if $T$ is a UGOA-tree, then $T^{\prime }$ is a UGOA-tree.
\end{itemize}
\end{ob}

\noindent\begin{proof}
$i)$ Since $y\in D$ and $v\in D\cap D^{\prime },$ it follows that $%
D-\left
\{ y\right \} $ is a GOA-set of $T^{\prime }$ and so $\left \vert
D^{\prime }\right \vert \leq \left \vert D\right \vert -1$. Moreover, since $v\in
D^{\prime }$, $D^{\prime }$ can be extended to a GOA-set of $T$ by adding $%
y$. Then $\left \vert D\right \vert \leq \left \vert D^{\prime }\cup
\{y\} \right
\vert =\left \vert D^{\prime }\right \vert +1$ and equality
holds.\vspace{2mm}\\
$ii)$ Since $D\cap V(T^{\prime })=D-\{y\}$ is a GOA-set of $T^{\prime }$
with cardinality $\left \vert D\right \vert -1=\left \vert D^{\prime
}\right
\vert ,$ $D\cap V(T^{\prime })$ is a $\gamma _{o}\left( T^{\prime
}\right) $-set.\vspace{2mm}\\
$iii)$ Let $B=\{y\}.$ In view of item $(i),$ $D^{\prime }$ can be extended
to a $\gamma _{o}\left( T\right) $-set by adding the unique vertex of $B.$
This and item $(ii)$ together with the uniqueness of $D$ imply that $%
D^{\prime }=D\cap V(T^{\prime })$ is the unique $\gamma _{o}$-set of $%
T^{\prime }.$
\end{proof}

\begin{ob}
\label{observation 4} Let $k$ be a positive integer and let $T$ be a tree
obtained from a nontrivial tree $T^{\prime }$ by adding $kP_{2}$  joining $k$ pairwise non-adjacent vertices of $kP_{2}$ to
the same leaf $v$ of $T^{\prime }.$ Let $w$ be the support vertex adjacent
to $v$, and let $D$ and $D^{\prime }$ be $\gamma _{o}\left( T\right) $-sets
of $T$ and $T^{\prime },$ respectively. If $w\in D\cap D^{\prime },$ then
the following three properties are satisfies.

\begin{itemize}
\item[$(i)$] $\left \vert D^{\prime }\right \vert =\left \vert D\right \vert
-k,$

\item[$(ii)$] $D\cap V(T^{\prime })$ is a $\gamma _{o}\left( T^{\prime
}\right) $-set,

\item[$(iii)$] if $T$ is a UGOA-tree, then $T^{\prime }$ is a UGOA-tree.
\end{itemize}
\end{ob}

\noindent\begin{proof}
Let $V(kP_{2})=\{x_{1},x_{2},\ldots,x_{k},y_{1},y_{2},\ldots,y_{k}\}$ and $%
E(kP_{2})=\{x_{i}y_{i}:i=1,2,\ldots,k\}.$ Let $v$ be a leaf of $T^{\prime }$
and $w$ be the support vertex adjacent to $v$. We assume that for each $%
i\in \{1,\ldots,k\},$ $y_{i}$ is adjacent to $v$ in $T.$\vspace{2mm}\\
$i)$ Obviously, all vertices of $\cup _{j=1}^{k}\left \{ y_{j}\right \} $
are support vertices in $T.$ Hence, in view of Observation \ref{observation
1} $(i)$, we can assume that $D$ contains all vertices of $\cup
_{j=1}^{k}\left \{ y_{j}\right \} .$ Therefore, since $w\in D,$ $D-(\cup
_{j=1}^{k}\left \{ y_{j}\right \} )$ is a GOA-set of $T^{\prime },$ which
means that $\left \vert D^{\prime }\right \vert \leq \left \vert D-(\cup
_{j=1}^{k}\left \{ y_{j}\right \} )\right \vert =\left \vert D\right \vert
-k.$ Observe that since $w\in D^{\prime },$ $D^{\prime }$ can be extended to
a GOA-set of $T$ by adding all vertices of $\cup _{j=1}^{k}\left \{
y_{j}\right \} .$ Hence $\left \vert D\right \vert \leq \left \vert
D^{\prime }\cup (\cup _{j=1}^{k}\left \{ y_{j}\right \} )\right \vert
=\left
\vert D^{\prime }\right \vert +k$ and so equality holds.\vspace{2mm}\\
$ii)$ The proof is similar to that of Observation \ref{observation 3}$(ii)$,
by taking $D\cap V(T^{\prime })=D-(\cup _{j=1}^{k}\left \{ y_{j}\right \} ).$\vspace{2mm}\\
$iii)$ The proof is similar to that of $(iii)$ of Observation \ref{observation 3}$(iii)$, by taking $%
B=\cup _{j=1}^{p}\{y_{j}\}.$
\end{proof}

\begin{ob}
\label{observation 5} Let $V(T^{\prime })$ be the vertex-set of a nontrivial
tree $T^{\prime },$ and let $D^{\prime }$ be a $\gamma _{o}(T^{\prime })$%
-set such $V(T^{\prime })-D^{^{\prime }}$ has a vertex $w$ with degree $%
q\geq 2$ and $\left \vert N_{T^{\prime }}(w)\cap (V(T^{\prime })-D^{\prime
})\right \vert \leq 1.$ Let $p$ be a positive integer such that
\begin{equation}
\left \{
\begin{array}{ll}
&p\leq q-1\text{ if }\left \vert N_{T^{\prime }}(w)\cap (V(T^{\prime
})-D^{\prime })\right \vert =0 , \\
\text{or} &\\
&p\leq q-3\text{ if }\left \vert N_{T^{\prime }}(w)\cap (V(T^{\prime
})-D^{\prime })\right \vert =1.%
\end{array}%
\right.   \label{cond1}
\end{equation}%
Let $T$ be a tree obtained from $T^{\prime }$ by adding $p$ subdivided stars
$SS_{k_{1}},\ldots,SS_{k_{p}}$ ($k_{i}\geq 2$ for all $i$) with centers $x_{1},$
$x_{2},\ldots,x_{p},$ respectively, and joining each $x_{i}$ $(1\leq i\leq p)$
at $w.$ Let $D$ be a $\gamma _{o}$-set of $T$. If $w$ and $%
x_{1},x_{2},\ldots,x_{p}$ are not in $D,$ then the following three properties
are satisfied.

\begin{itemize}
\item[$(i)$] $\left \vert D^{\prime }\right \vert =\left \vert D\right \vert
-\overset{p}{\underset{i=1}{\sum }}k_{i},$

\item[$(ii)$] $D\cap V(T^{\prime })$ is a $\gamma _{o}\left( T^{\prime
}\right) $-set,

\item[$(iii)$] if $T$ is a UGOA-tree, then $T^{\prime }$ is also a UGOA-tree.
\end{itemize}
\end{ob}

\noindent\begin{proof}
For $i\in \{1,\ldots,p\},\ $let $S\left( SS_{k_{i}}\right) $ be a support
vertex-set of $SS_{k_{i}}.$\vspace{2mm}\\
$i)$ Since $w$ together with $x_{1},x_{2},\ldots,x_{p}$ are not in $D,$ all
vertices of $\cup _{i=1}^{p}S\left( SS_{k_{i}}\right) $ must be in $D.$\
Therefore,  $D\backslash \overset{p}{\underset{i=1}{\cup }}S\left(
SS_{k_{i}}\right) $ is a GOA-set of $T^{\prime },$ giving that $\left \vert
D^{\prime }\right \vert \leq \left \vert D\right \vert -\overset{p}{\underset%
{i=1}{\sum }}k_{i}.$\\ On the other hand, let $A=\cup _{i=1}^{p}S\left( SS_{k_{i}}\right) \cup
D^{\prime }.$ We have to show that $A$ is a GOA-set of $T.$ For this, it
suffices to show that $\left \vert N_{T}\left[ z\right] \cap A\right \vert
\geq \left \vert N_{T}\left[ z\right] -A\right \vert $ for each $z\in
\{w,x_{1},x_{2},\ldots,x_{p}\}.$ Indeed,  we have to distinguish between two cases.\\
\noindent \textbf{Case 1.} $z=x_{i},$ for some $i\in \{1,\ldots,p\}.$ \\We have then
\begin{equation*}
\left \vert N_{T}\left[ z\right] \cap A\right \vert =\left \vert N_{T}\left[
z\right] \cap \cup _{i=1}^{p}S\left( SS_{k_{i}}\right) \right \vert
=k_{i}\geq 2,
\end{equation*}%
and
\begin{equation*}
\left \vert N_{T}\left[ z\right] -A\right \vert =\left \vert \{z,w\} \right
\vert =2.
\end{equation*}

\noindent\textbf{Case 2.} $z=w.$ \\We have then
\begin{equation*}
\left \vert N_{T}\left[ z\right] \cap A\right \vert =\left \{
\begin{array}{ccl}
q&\text{if}&\left \vert N_{T^{\prime }}(w)\cap (V(T^{\prime })-D^{\prime
})\right \vert =0, \\
q-1&\text{if}&\left \vert N_{T^{\prime }}(w)\cap (V(T^{\prime })-D^{\prime
})\right \vert =1.%
\end{array}%
\right.
\end{equation*}%
and
\begin{equation*}
\left \vert N_{T}\left[ z\right] -A\right \vert =\left \{
\begin{array}{ccl}
p+1&\text{if}&\left \vert N_{T^{\prime }}(w)\cap (V(T^{\prime })-D^{\prime
})\right \vert =0, \\
p+2&\text{if}&\left \vert N_{T^{\prime }}(w)\cap (V(T^{\prime })-D^{\prime
})\right \vert =1.%
\end{array}%
\right.
\end{equation*}

According to (\ref{cond1}), we have in each case $\left \vert N_{T}\left[ z%
\right] \cap A\right \vert \geq \left \vert N_{T}\left[ z\right] -A\right \vert
$ for each $z\in \{w,x_{1},x_{2},\ldots,x_{p}\}.$ Therefore $A$ is a GOA-set of
$T,$ giving that $\left \vert D\right \vert \leq \left \vert A\right \vert
=\left \vert D^{\prime }\right \vert +\overset{p}{\underset{i=1}{\sum }}k_{i}$. Hence the equality holds.\vspace{2mm}\\
$ii)$ Using the fact that $D\cap V(T^{\prime })=D\backslash\cup
_{i=1}^{p}S\left( SS_{k_{i}}\right),$ this property follows in a similar manner as the proof of Observation \ref{observation 3}$(ii)$.\vspace{2mm}\\
$(iii)$ This property follows in a similar manner as the proof of Observation \ref%
{observation 3}$(iii)$, by taking $B=\cup _{i=1}^{p}S\left(
SS_{k_{i}}\right) .$
\end{proof}

\section{The main result}

In order to characterize the trees with unique minimum global offensive
alliance, we define a family $\mathcal{F}$ of all trees $T$ that can be
obtained from a sequence $T_{1},T_{2},\ldots,T_{r}$ $\left( r\geq 1\right) $ of
trees, where $T_{1}$ is the path $P_{3}$ centered at a vertex $y,$ $T=T_{r},$
and if $r\geq 2,$ $T_{i+1}$ is obtained recursively fom $T_{i}$ by one of
the following operations. Let $A\left( T_{1}\right) =\left \{ y\right \} .$

\begin{itemize}
\item Operation $\mathcal{O}_{1}:$ Attach a vertex by joining it to any
support vertex of $T_{i}.$ Let $A\left( T_{i+1}\right) =A\left( T_{i}\right).$
\item Operation $\mathcal{O}_{2}:$ Attach a path $P_{3}=u$-$v$-$w$ by
joining $v$ to any support vertex of $T_{i}.$ Let $A\left( T_{i+1}\right)
=A\left( T_{i}\right) \cup \left \{ v\right \} .$
\item Operation $\mathcal{O}_{3}:$ Let $w$ be a support vertex of $T_{i}$
that satisfies one of the following two conditions.
\begin{itemize}
\item[1.] $l_{T_{i}}(w)\geq 3,$
\item[2.] $\left \vert N_{T_{i}}[w]\cap A(T_{i})\right \vert <\left \vert
N_{T_{i}}(w)\cap (V(T_{i})-A(T_{i})\right \vert $ or
\begin{itemize}
\item either $l_{T_{i}}(w)=2$ and $N_{T_{i}}(w)-A(T_{i})$ has a vertex $%
w_{t} $ such that $\left \vert N_{T_{i}}(w_{t})\cap A(T_{i})\right \vert
\leq \left \vert N_{T_{i}}[w_{t}]\cap (V(T_{i})-A(T_{i})\right \vert +1,$
\item or $l_{T_{i}}(w)=1$ and $N_{T_{i}}(w)-A(T_{i})$ has two vertices $%
w_{p},w_{q}$ so that for $l=p,q,$ $\left \vert N_{T_{i}}(w_{l})\cap
A(T_{i})\right \vert \leq \left \vert N_{T_{i}}[w_{l}]\cap
(V(T_{i})-A(T_{i})\right \vert +1.$
\end{itemize}
\end{itemize}
\end{itemize}

Let $kP_{2}$ be the disjoint union of $k\geq 1$ copies of $P_{2}$, and let $%
B $ be a set of $k$ pairwise non-adjacent vertices of $kP_{2}.$ Add $kP_{2}$
and attach all vertices of $B$ to a same leaf in $T_{i}$ that is adjacent to
$w.$ Let $A\left( T_{i+1}\right) =A\left( T_{i}\right) \cup B.$

\begin{itemize}
\item Operation $\mathcal{O}_{4}:$ Let $w\in V\left( T_{i}\right) -A\left(
T_{i}\right) $ be a vertex of degree $q\geq 2$ in $T_{i}$ such that $%
\left
\vert N_{T_{i}}(w)\cap (V\left( T_{i}\right) -A\left( T_{i}\right)
)\right
\vert \leq 1$. Attach $p\geq 1$ subdivided stars $SS_{k_{i}}$ ($k_{i}\geq 2$ for $1\leq i\leq p)$ with support vertex-set $S\left(
SS_{k_{i}}\right) $ and of center $x_{i}$ by joining $x_{i}$ to $w$ for all $%
i$ such that
\[p \leq \left\{\begin{array}{cl}
                 q -1 & \text{if } \left \vert N_{T_{i}}(w)\cap (V\left(T_{i}\right) -A\left( T_{i}\right) )\right \vert =0,\\
                 q-3 & \text{if } \left \vert N_{T_{i}}(w)\cap (V\left(T_{i}\right) -A\left( T_{i}\right) )\right \vert =1.
                    \end{array}
              \right.\]

\end{itemize}

Let $A\left( T_{i+1}\right) =A\left( T_{i}\right) \cup (\cup
_{i=1}^{p}S\left( SS_{k_{i}}\right) ).\bigskip $

Before stating our main result, we need the following lemma.

\begin{lemma}
\label{lemma 1}If $T\in {\cal F},$ then $A\left( T\right) $ is the
unique $\gamma _{o}\left( T\right) $-set.
\end{lemma}
\noindent\begin{proof}
Let $T\in {\cal F}$. We proceed by induction on the number of operations, say $r$, required to construct $T.$ The
property is true if $T$ is a path $P_{3}$ centered at $y$ since $A\left(
T\right) =\left \{ y\right \} $ is the unique $\gamma _{o}\left( T\right)$-set. This establishes the base case. \\
Assume that for any tree $T^{\prime }\in {\cal F} $ that can be
constructed with $r-1$ operations, $A\left( T^{\prime }\right) $ is the
unique $\gamma _{o}\left( T^{\prime }\right) $-set. Let $T=T_{r}$ with $%
r\geq 2$ and $T^{\prime }=T_{r-1}.$ We distinguish between four cases.\vspace{2mm}\\
\noindent{\bf Case 1.} $T$ is obtained from $T^{\prime }$ by using Operation
$\mathcal{O}_{1}$.\\ Assume that $T$ is obtained from $T^{\prime }$ by
attaching an extra vertex at a support vertex $u$ of $T^{\prime }$. In view of Observation \ref{observation 1} $(ii),$ $u\in A(T^{\prime })$. Hence $A(T^{\prime })$ can be extended to a GOA-set of $T$. By Observation %
\ref{observation 2} $(i),$ $\gamma _{o}\left( T\right) =\gamma _{o}\left(
T^{\prime }\right) ,$ implying that $A(T^{\prime })$ is a $\gamma _{o}\left(
T\right) $-set. Applying the inductive hypothesis to $T^{\prime },$ $%
A(T^{\prime })$ is the unique $\gamma _{o}\left( T^{\prime }\right) $-set.
It follows that $A\left( T\right) =A\left( T^{\prime }\right) $ is the
unique $\gamma _{o}\left( T\right)$-set.\vspace{2mm}\\
\noindent{\bf Case 2.} $T$ is obtained from $T^{\prime }$ by using Operation
$\mathcal{O}_{2}$. \\  $A\left( T^{\prime }\right) \cup \left \{
v\right \} $ is a GOA-set of $T$. By Observation \ref{observation 3} $(i),$
$\gamma _{o}\left( T\right) =\gamma _{o}\left( T^{\prime }\right) +1,$
meaning that $A\left( T^{\prime }\right) \cup \left \{ v\right \} $ is a $%
\gamma _{o}\left( T\right) $-set. The inductive hypothesis sets that $%
A\left( T^{\prime }\right) $ is the unique $\gamma _{o}\left( T^{\prime
}\right) $-set. Thus $A\left( T\right) =A\left( T^{\prime }\right) \cup
\left \{ v\right \} $ is the unique $\gamma _{o}\left( T\right)$-set. \vspace{2mm}\\
\noindent {\bf Case 3.} $T$ is obtained from $T^{\prime }$ by using Operation
$\mathcal{O}_{3}$. \\  $A\left( T^{\prime }\right) \cup B$ is a GOA-set
of $T$. Observation \ref{observation 4} $(i)$ sets that $\gamma _{o}\left(
T\right) =\gamma _{o}\left( T^{\prime }\right) +k,$ which means that $%
A\left( T^{\prime }\right) \cup B$ is a $\gamma _{o}\left( T\right) $-set.
By the inductive hypothesis, $A\left( T^{\prime }\right) $ is the unique $%
\gamma _{o}\left( T^{\prime }\right) $-set. Thus $A\left( T\right) =A\left(
T^{\prime }\right) \cup B$ is the unique $\gamma _{o}\left( T\right)$-set.\vspace{2mm}\\
\noindent{\bf Case 4.} $T$ is obtained from $T^{\prime }$ by using Operation
$\mathcal{O}_{4}$. \\ $A\left( T^{\prime }\right) \cup (\cup
_{i=1}^{p}S\left( SS_{k_{i}}\right) )$ is a GOA-set of $T$. According to
Observation \ref{observation 5} $(i)$, we have $\gamma _{o}\left( T\right)
=\gamma _{o}\left( T^{\prime }\right) +\sum\nolimits_{i=1}^{p}k_{i},$
whence, $A\left( T^{\prime }\right) \cup (\cup _{i=1}^{p}S\left(
SS_{k_{i}}\right) )$ is a $\gamma _{o}\left( T\right) $-set. By the
inductive hypothesis, $A\left( T^{\prime }\right) $ is the unique $\gamma
_{o}\left( T^{\prime }\right) $-set. It follows that $A\left( T\right)
=A\left( T^{\prime }\right) \cup (\cup _{i=1}^{p}S\left( SS_{k_{i}}\right) )$
is the unique $\gamma _{o}\left( T\right) $-set.
\end{proof}
\\

Remark that in each case, $A(T_{i+1})$ is obtained from $A(T_{i})$ by adding
all support vertices in $T_{i+1}\backslash T_{i}.$ Hence the following
corollary is immediate.

\begin{corollary}
Let $T\in {\cal F} $ and $S(T)$ be a set of support vertices in $T$. Then
$\gamma_{o}\left( T\right) \geqslant \left \vert S(T)\right \vert.$
\end{corollary}

Now we are ready to prove our main result.

\begin{theorem}
\label{thm}A tree $T$ is a UGOA-tree if and only if $T=K_{1}$ or $T\in
{\cal F}.$
\end{theorem}

\noindent\begin{proof}
It is obvious that $T=K_{1}$ is a UGOA-tree. Also, Lemma \ref{lemma 1}
states that any member of ${\cal F} $ is a UGOA-tree. Now, we prove the
converse by induction on the number $n$ of vertices of $T$. The converse
holds trivially for $n=1$ and $3$ but not for $n=2$ since $P_{2}$ is not a
UGOA-tree. When $n=4,$ $T$ is either a $K_{1,3}$ or a $P_{4}$. Clearly $P_{4}
$ is not a UGOA-tree, whilst $K_{1,3}$ is a UGOA-tree that can be obtained
from a $P_{3}$ using operation $\mathcal{O}_{1}$, and so $K_{1,3}\in
{\cal F} $. If $n=5$, then $T$ is either a double star $S_{1,2}$ which is
not a UGOA-tree, or it is a $K_{1,4}$ or $P_{5}$ that are UGOA-tree since $%
K_{1,4}$ can be obtained from $K_{1,3}$ by using operation $\mathcal{O}_{1},$
and $P_{5}$ can be obtained from a $P_{3}$ by using operation $\mathcal{O}%
_{3}.$ Therefore $K_{1,4}$ and $P_{5}$ are in ${\cal F} $. This
establishes the base case.

Now, let $n\geq 6$ and assume that any tree $T^{\prime }$ of order $3\leq
n^{\prime }<n$ with the unique $\gamma _{o}\left( T^{\prime }\right) $-set
is in ${\cal F} $. Let $T$ be a tree of order $n$ with the unique $\gamma
_{o}\left( T\right) $-set $D$ and let $s\in S(T).$ By Observation \ref%
{observation 1} $(ii),$ $s\in D.$\ If $l_{T}\left( s\right) \geq 3,$ then
let $T^{\prime }$ be the tree obtained from $T$ by removing a leaf adjacent
to $s$ and let $D^{\prime }$ be a $\gamma _{o}\left( T^{\prime }\right) $%
-set. Then, clearly $n^{\prime }=\left \vert V(T^{\prime })\right \vert
=n-1\geq 5,$ and $l_{T^{\prime }}\left( s\right) \geq 2,$ so $s\in D^{\prime
}$ by Observation \ref{observation 1} $(iii).$ According to Observation \ref%
{observation 2} $(ii)$, $T^{\prime }$ is UGOA-tree. Applying the inductive
hypothesis to $T^{\prime },$ we get $T^{\prime }\in {\cal F} .$ Thus $T$
is obtained from $T^{\prime }$ by operation $\mathcal{O}_{1},$ implying that
$T\in {\cal F} $. Assume now that
\begin{equation}
\text{for each }x\in S(T)\text{, }l_{T}(x)\leq 2.  \label{xx0}
\end{equation}

Root $T$ at a vertex $r$ of maximum eccentricity.\ Let $u$ be a support
vertex of maximum distance from $r\ $and let $u^{\prime }$ be a leaf
adjacent to $u$. Let $v$ and $w$ be the parents of $u$ and $v,$
respectively, in the rooted tree. We consider two cases. \vspace{2mm}\\
\noindent{\bf Case 1.} $v\in D.$ \\If $l_{T}\left( u\right) =1,$ then $D\cup \left \{ u^{\prime }\right \}
-\left \{ u\right \} $ is a $\gamma _{o}(T)$-set, contradicting the uniqueness
of $D$ as a $\gamma _{o}(T)$-set. Hence by (\ref{xx0}), $l_{T}\left(
u\right) =2.$ We claim that $v\in S(T).$ Suppose not. Then either $w\in D$
and so $D-\left \{ v\right \} $ is a GOA-set of $T$ with cardinality less than
$\left \vert D\right \vert ,$ contradicting the minimality of $D,$\ or $%
w\notin D$ and so $D-\left \{ v\right \} \cup \left \{ w\right \} $ is a $\gamma
_{o}\left( T\right) $-set, contradicting the uniqueness of $D$ as a $\gamma
_{o}(T)$-set. This completes the proof of the claim. Let $T^{\prime }=T-T_{u}
$ and $D^{\prime }$ be a $\gamma _{o}$-set of $T^{\prime }.$ By Observation %
\ref{observation 1}$(i)$, we can assume that $D^{\prime }$ contains all
support vertices in $T^{\prime }.$ Since $\left \vert V(T_{u})\right \vert =3,$
it follows that $n^{\prime }=\left \vert V(T^{\prime })\right \vert =n-3\geq 3$
and so $T^{\prime }\neq P_{2}.$ By Observation \ref{observation 3}$(iii)$, $%
T^{\prime }$ is a UGOA-tree. Applying our inductive hypothesis, we get $%
T^{\prime }\in {\cal F} $. Thus, $T$ can be obtained from $T^{\prime }$
by operation $\mathcal{O}_{2}$ and so $T\in {\cal F}.$ \vspace{2mm}\\
\noindent{\bf Case 2.}  $v\notin D.$ \\According to Observation \ref{observation 1}$(ii)$, $v\notin S(T)$ and so $%
l_{T}(v)=0.$ Let $k=\left \vert N_{T}(v)-\{w\} \right \vert.$ We have then $%
d_{T}(v)=k+1$ and  since $u\in N_{T}(v)-\{w\}$, we clearly deduce $k\geq 1$. For $i\in
\{1,\ldots,k\},$ let $u_{i}\in N_{T}(v)-\{w\}$ such that $u_{1}=u$. The choice
of $v$ sets that
\begin{equation}
u_{i}\in S(T),\text{ }l_{T}(u_{i})\geq 1\text{ and so }u_{i}\in D\text{ for
all }i.  \label{xx00}
\end{equation}%
Hence by (\ref{xx0}), we have $1\leq l_{T}(u_{i})\leq 2$ for all $i.$ Assume
first that $l_{T}(u_{j})=2$ for some $j$ in $\{1,\ldots,k\}.$ Without loss of
generality, let $j=1.$ Then $u$ has a further neighbor $u^{\prime \prime
}\neq u^{\prime }$ in $T.$ Let $T^{\prime }=T-\{u^{\prime \prime }\}$ and $%
D^{\prime }$ be any $\gamma _{o}$-set of $T^{\prime }.$ Clearly $u^{\prime }$
is the unique leaf of $u$ in $T^{\prime }.$ We claim that $u\in D^{\prime }.$
Suppose not. Then $u^{\prime }$ and $v$ must be in $D^{\prime }$ and
therefore $D^{\prime \prime }=(D^{\prime } \backslash\{u^{\prime }\})\cup
\{u\}$ is a further $\gamma _{o}\left( T\right) $-set other than $D$ (since $%
v$ belongs to $D^{\prime \prime }$ and not to $D),$ a contradiction. This
completes the proof of the claim. We have $n^{\prime }=n-1\geq 5.$ By
Observation \ref{observation 2}$(iii),$ $T^{\prime }$ is a UGOA-tree.
Applying our inductive hypothesis to $T^{\prime },$ we get $T^{\prime }\in
{\cal F} .$ Hence $T$ is obtained from $T^{\prime }$ by operation $%
\mathcal{O}_{1}$, implying that$\ T\in {\cal F} .$ Assume now that
\begin{equation}
l_{T}(u_{i})=1\text{ and hence }d_{T}(u_{i})=2\text{ for all }i.  \label{xx1}
\end{equation}%
For all $i\in \{1,\ldots,k\}$, let $u_{i}^{\prime }$ be the unique leaf
adjacent to $u_{i}$ (with $u_{1}^{\prime }=u^{\prime }$). We distinguish
between two subcases, depending on whether $w$ belongs to $D$ or not.\vspace{2mm}\\
\noindent{\bf Case 2.1.} $w\in D.$ \\In view of (\ref{xx1}), $T_{v}-\left \{ v\right \} =kP_{2}$ with $V(kP_{2})=%
\{u_{1},u_{2},\ldots,u_{k},u_{1}^{\prime },u_{2}^{\prime },\ldots,u_{k}^{\prime
}\} $ and $E(kP_{2})=\{u_{i}u_{i}^{\prime }:i=1,2,\ldots,k\}.$ Let $T^{\prime
}=T-\left( T_{v}-\left \{ v\right \} \right) .$ Clearly $v\in L(T^{\prime })$
and $w\in S(T^{\prime }).$ If $n^{\prime }=\left \vert V(T^{\prime
})\right
\vert =2,$ then $T$ is a wounded spider with exactly one
non-subdivided edge and in this case, it is not difficult to see that such a graph
is not a UGOA-tree. Hence assume that $n^{\prime }\geq 3.$ We claim the
following:\\
If $l_{T}(w)\in \{0,1\},$ then one of the two conditions
holds:
\begin{itemize}
\item[C$_{1}:$] $\left \vert N_{T}[w]\cap D\right \vert \leq \left \vert
N_{T}(w)\cap (V(T)-D\right \vert.$
\item[C$_{2}:$] $(i)$ either $l_{T}(w)=1$ and $N_{T}(w)-D$ has a vertex $%
w_{t}$ such that \[\left \vert N_{T}(w_{t})\cap D\right \vert \leq \left \vert
N_{T}[w_{t}]\cap (V(T)-D)\right \vert +1 \]\\
$(ii)$ or, $l_{T}(w)=0$ and $N_{T}(w)-D$ has two vertices $w_{p},w_{q}$ such
that for $l \in \{p, q\}$,  \[\left \vert N_{T}(w_{l})\cap D\right \vert \leq \left \vert
N_{T}[w_{l}]\cap (V(T)-D)\right \vert +1.\]
\end{itemize}
Indeed, suppose that $C_{1}$ and $C_{2}$ are not satisfied. Assume
first that $l_{T}(w)=1,$ so $L_{T}(w)$ has exactly one vertex, say $%
w^{\prime }$. In this case $D-\left \{ w\right \} \cup \left \{ w^{\prime
}\right \} $ is a $\gamma _{o}\left( T\right) $-set different from $D$, a
contradiction. Now, assume that $l_{T}(w)=0.$ Since $C_{2}$ is not
fulfilled, item $(ii)$ of $C_{2}$ is satisfied for at most one vertex in $%
N_{T}(w)-D$, say  $w^{\prime \prime }$. Then $%
D-\left \{ w\right \} \cup \left \{ w^{\prime \prime }\right \} $ is a $%
\gamma _{o}\left( T\right) $-set different from $D$, a contradiction. If no
vertex in $N_{T}(w)-D$ for which item $(ii)$ of $C_{2}$ is satisfied, then $%
D-\left \{ w\right \} \cup \left \{ v\right \} $ is a $\gamma _{o}\left(
T\right) $-set different from $D$, which leads to a contradiction again.
This complete the proof of the claim.\\
Observe that when $l_{T^{\prime }}(w)\in \{1,2\},$  the previous claim remain true by
replacing $D$ by $D^{\prime }$ and $T$ by $T^{\prime }.$ Thus, according to
Observation \ref{observation 4} $(iii), T^{\prime }$ is a UGOA-tree.
By induction on $T^{\prime },$ we get $T^{\prime }\in {\cal F} .$ Since $T$ is obtained from $T^{\prime }$ by using operation $%
\mathcal{O}_{3},$ we directly obtain $T\in
{\cal F}$.  \vspace{2mm}\\
\noindent\textbf{Case 2.2.} $w$ $\notin D.$\\
By Observation \ref{observation 1}$(ii)$, $w\notin S(T)$ and so $%
l_{T}(w)=0. $ Since $v$ and $w$ are in $V(T)-D,$ $v$ must have at least two
neighbors in $D$. Hence $d_{T}(v)=k+1\geq 3.$ Let $t$ be the parent of $w,$
and let $X,Y$ and $Z$ be the following sets
\begin{equation*}
Y=C(w)\cap S(T),\text{ }X=C(w)-Y\text{ and }Z=D(w)\cap \left( S(T)-Y\right) .
\end{equation*}%
Observe that $v\in X,$ $u\in Z,$ $N_{T}(w)=\{t\} \cup X\cup Y$ and every
vertex in $Z$ plays the same role as $u.$ Therefore by (\ref{xx00}), we have
$Z\subset D$ since $Z\subset S(T),$ and by (\ref{xx1}), every vertex in $Z$
has exactly two neighbors such that one of them is a leaf and the other one
is in $X$. Furthermore, as $v\in X,$ $u_{i}\in Z$ for all $i\in \{1,\ldots,k\},$
so $\left \vert Z\right \vert \geq k\geq 2.$ Notice also that $\left \vert
X\right \vert \geq 1$ since $v\in X.$ Likewise $\left \vert Y\right \vert
\geq 1 $ since $D$ is a $\gamma _{o}(T)$-set. \ It is clear that $Y\subseteq
S(T)\ $and thus $Y\subseteq D$ by Observation \ref{observation 1}$(ii).$
Setting
\begin{equation*}
X=\{x_{1},x_{2},\ldots x_{p}\}(p\geq 1)\text{ with }x_{1}=v\text{ and }\left
\vert Y\right \vert =q-1\text{ }(q\geq 2).
\end{equation*}%
Since every vertex in $X$ plays the same role as $v,$ $x_{i}\in V(T)-D$ for
all $i\in \{1,\ldots,p\}.$  Setting%
\begin{equation*}
p_{i}=\left \vert N_{T}(x_{i})-\{w\} \right \vert \text{ for }i=1,\ldots,p.
\end{equation*}%
Then $p_{1}=k.$ Since for all $i \in \{1,\ldots,p\}$, $x_{i}$  and $w$ are in $V(T)-D,$
$x_{i}$ must have at least two neighbors in $Z$. Hence $d_{T}(x_{i})=p_{i}+1%
\geq 3$. This means that for all $i\in
\{1,\ldots,p\},$ $V(T_{x_{i}})$ induces a subdivided star $SS_{p_{i}}$ of order
$p_{i}+1$ centered at $x_{i}$. Since $w\in V(T)-D,$ inequality (\ref{def})
is valid by replacing $v$ with $w.$ This gives
\begin{equation}
p\leq q-1\text{ if }t\in D,\text{ or }p\leq q-3\text{ otherwise.}  \label{pq}
\end{equation}%
Let $T^{\prime }=T-\cup (\cup _{i=1}^{p}T_{x_{i}})$ and $D^{\prime }$ be a $%
\gamma _{o}(T^{\prime })$-set. Observe that $T^{\prime }$ contains at least
one $P_{3}$ as an induced subgraph, which means that $n^{\prime
}=\left
\vert V(T^{\prime })\right \vert \geq 3.$ For all $i\in \{1,\ldots,p\}$%
, let $S(SS_{p_{i}})$ be the support vertex-set of $SS_{p_{i}}.$ Clearly $%
\cup _{i=1}^{p}S(SS_{p_{i}})=Z$ and $N_{T^{\prime }}(w)=Y\cup \{t\},$ so
\begin{equation*}
d_{T^{\prime }}(w)=q\geq 2.
\end{equation*}%
According to Observation \ref{observation 1} $(i)$, we can assume that $%
Y\subset D^{\prime }$ since $Y\subset S(T^{\prime }).$ Then $t$ is the only
neighbor of $w$ in $T^{\prime }$ that may not be in $D^{\prime }$, that is%
\begin{equation*}
\left \vert N_{T^{\prime }}(w)\cap (V(T^{\prime })-D^{\prime })\right \vert
\leq 1.
\end{equation*}%
If $t\in D^{\prime },$ then the minimality of $D^{\prime }$ sets that $w\in
V(T^{\prime })-D^{\prime },$ because otherwise,  we replace $w$
by $t$ in $D^{\prime }$. \\
By Observation \ref{observation 5} $(ii)$ and $(iii),$ we have $D^{\prime
}=D\cap V(T^{\prime }).$ Hence $t\in D$ if and only if $t\in D^{\prime }$.
Notice that if $t\in D^{\prime },$ then $N_{T^{\prime }}(w)\cap (V(T^{\prime
})-D^{\prime })$ is an empty-set, otherwise, $t$ would be the unique vertex
of $N_{T^{\prime }}(w)\cap (V(T^{\prime })-D^{\prime }).$ Thus (\ref{pq})
can be rewritten as follows.
\begin{equation*}
\text{If }\left \vert N_{T^{\prime }}(w)\cap (V(T^{\prime })-D^{\prime
})\right \vert =0,\text{ then }p\leq q-1,\text{ }
\end{equation*}%
and
\begin{equation*}
\text{if }\left \vert N_{T^{\prime }}(w)\cap (V(T^{\prime })-D^{\prime
})\right \vert =1,\text{ then }p\leq q-3.
\end{equation*}

Again Observation \ref{observation 5}$(iii)$ sets that $T^{\prime }$ is a
UGOA-tree. Applying the inductive hypothesis to $T^{\prime },$ we deduce $%
T^{\prime }\in {\cal F} .$ Now since $T$ can be obtained from $T^{\prime
} $ by operation $\mathcal{O}_{4}$, and finally $T\in {\cal F} .$ This completes the
proof of Theorem \ref{thm}.\bigskip
\end{proof}

\section{Open Problems}
The previous results motivate the following problems.

\begin{itemize}
\item[1-] Characterize other UGOA-graphs.

\item[2-] Characterize trees with unique minimum defensive alliance sets
(UGDA).
\end{itemize}

\end{document}